\documentclass[11pt,centertags,reqno]{amsart}
\usepackage{amsfonts}
\usepackage{amssymb}
\usepackage{graphicx}
\usepackage{amscd}

\theoremstyle{definition}

\theoremstyle{plain}
\numberwithin{equation}{section}
\input{tcilatex}

\begin{document}

\begin{center}
{\LARGE Rates of convergence for a class of generalized quasi contractive
mappings in Kohlenbach hyperbolic spaces\bigskip }

Zahid Akhtar$^{a}$ and Muhammad Aqeel Ahmad Khan$^{b,\ast }$

$^{a}$Department of Mathematics, Govt. S. E. College Bahawalpur,

Bahawalpur, 63100, Pakistan

$^{b}$Department of Mathematics, COMSATS Institute of Information Technology
Lahore,

Lahore, 54000, Pakistan

February 23, 2018 \let\thefootnote\relax\footnote{%
*Corresponding author\newline
E-mail addresses: (Z. Akhter) {\small zahid\_9896@yahoo.com, }(M.A.A. Khan)
itsakb@hotmail.com, maqeelkhan@ciitlahore.edu.pk}
\end{center}

\textbf{Abstract}: This paper is a continuation to the study of generalized
quasi contractive operators, essentially due to Akhtar et al. [A multi-step
implicit iterative process for common fixed points of generalized $C^{q}$%
-operators in convex metric spaces, Sci. Int., 25(4) (2013), 887-891], in
spaces of nonpositive sectional curvature. We aim to establish results
concerning convergence characteristics of the classical iterative algorithms
such as Picard, Mann, Ishikawa and Xu-Noor iterative algorithms associated
with the proposed class of generalized quasi contractive operators.
Moreover, we adopt the concept introduced by Berinde [Comparing
Krasnosel'skii and Mann iterative methods for Lipschitzian generalized
pseudo-contractions, Int. Conference on Fixed Point Theory Appl., 15-26,
Yokohama Publ., Yokohama, 2004.] for a comparison of the corresponding rates
of convergence of these iterative algorithms in such setting of spaces. The
results presented in this paper improve and extend some recent corresponding
results in the literature.\\[1mm]

\section{Introduction}

Fixed point theory (FPT) contributes significantly to the theory of
nonlinear functional analysis. Iterative algorithms, with respect to various
nonlinear mappings, are ubiquitous in FPT and have been successfully applied
in the study of a variety of nonlinear phenomena. The theory of iterative
construction of fixed points of a nonlinear mapping under suitable set of
control conditions is coined as metric fixed point theory (MFPT). MFPT is a
fascinating field of research and has emerged as a powerful tool to solve
various nonlinear real world problems, such as Fredholm and Volterra
integral equations, ordinary differential equations, partial differential
equations and image processing. MFPT has its roots in the celebrated Banach
Contraction Principle (BCP) which not only guarantees the existence of a
unique fixed point of a contraction but also describes an approximant for
the construction of such a unique fixed point. It is worth mentioning that
the BCP also gives a geometric rate of convergence for the classical Picard
iterative algorithm to the unique fixed point. The BCP is a frequently cited
result in the whole theory of analysis and dominates FPT for the class of
contractions.\medskip

It is worth mentioning that the simplicity and applicability of the BCP
paved the way for developing a new class of mappings satisfying generalized
contractive condition. Most of the generalizations of the BCP possess the
same characteristics regarding the existence of a unique fixed point which
can be constructed by the Picard iterative algorithm. However, there are
certain contractive type mappings for which the construction of fixed points
is also possible via Krasnosel'skii \cite{Krasnosel'skii}, Mann \cite%
{Groetsch,Mann}, Ishikawa \cite{Ishikawa} and Xu-Noor \cite{Xu-Noor}
iterative algorithms. In MFPT, different iterative algorithms can be
evaluated with respect to various characteristics, inter alia, convergence
characteristics and rates of convergence. The later concept has its own
importance in MFPT and therefore we adopt the concept introduced by Berinde 
\cite{Berinde1} for a comparison of the rates of convergence of different
iterative algorithms involving a nonlinear mapping.\medskip

Since a variety of problems corresponding to the real world nonlinear
phenomena can be transformed into fixed point problems (FPP). Therefore, it
is natural to study FPP associated with a class of mappings in a suitable
nonlinear framework. The term nonlinear framework for FPT is referred as a
metric space embedded with a "convex structure". It is remarked that the
non-positively curved hyperbolic space, introduced by Kohlenbach \cite%
{Kohlenbach}, provides rich geometrical structures suitable for MFPT of
various classes of mappings. For the results concerning MFPT in Kohlenbach
hyperbolic spaces, see, for example, \cite{Fukhar
FPTA,Fukhar-Zahid,KKF-NA,KFK FPTA,Khan JIA,Khan Fukhar FPTA,Khan-Fukhar-Amna}
and the references cited therein. We are, therefore, interested into
iterative construction of fixed points of the class of quasi contractive
mappings in Kohlenbach hyperbolic spaces. As a consequence, we establish
results concerning rates of convergence associated with the modified Mann,
Ishikawa and Xu-Noor iterative algorithms, involving the class of quasi
contractive mappings, in comparison to the classical Picard iterative
algorithm in Kohlenbach hyperbolic spaces.

\section{Preliminaries}

Throughout this paper, we work in the setting of hyperbolic spaces
introduced by Kohlenbach \cite{Kohlenbach} and hence the term Kohlenbach
hyperbolic spaces as one can find different notions of hyperbolic spaces in
the current literature, see \cite{Goebel Kirk,Goebel Reich,Reich
Zasalvski,Reich Shafrir}.

A Kohlenbach hyperbolic space $X$ is a metric space $(X,d)$ together with a
convexity mapping $W:$ $X^{2}\times \lbrack 0,1]\rightarrow X$ satisfying%
\begin{equation*}
\begin{array}{l}
\text{(W1) }d(u,W(x,y,\alpha ))\leq \alpha d(u,x)+(1-\alpha )d(u,y) \\ 
\text{(W2) }d(W(x,y,\alpha ),W(x,y,\beta ))=\left\vert \alpha -\beta
\right\vert d(x,y) \\ 
\text{(W3) }W(x,y,\alpha )=W(y,x,(1-\alpha )) \\ 
\text{(W4) }d(W(x,z,\alpha ),W(y,w,\alpha ))\leq \alpha d(x,y)+\left(
1-\alpha \right) d(z,w)\newline
\end{array}%
\end{equation*}%
for all $x,y,z,w\in X$ and $\alpha ,\beta \in \lbrack 0,1].$ A subset $K$ of
a hyperbolic space $X$ is convex if $W(x,y,\alpha )\in K$ for all $x,y\in K$
and $\alpha \in \lbrack 0,1].$ A hyperbolic space $X$ is uniformly convex 
\cite{L-JMAA} if for all $u,x,y\in X,$ $r>0$ and $\epsilon \in (0,2],$ there
exists $\delta \in (0,1]$ such that%
\begin{equation*}
d\left( W(x,y,\frac{1}{2}),u\right) \leq (1-\delta )r
\end{equation*}%
whenever $d(x,u)\leq r,d(y,u)\leq r$ and $d(x,y)\geq r\epsilon .$

A mapping $\eta :(0,\infty )\times (0,2]\rightarrow (0,1]$ providing such $%
\delta =\eta (r,\epsilon )$ for given $r>0$ and $\epsilon \in (0,2]$\ is
called modulus of uniform convexity. For more on hyperbolic spaces, we refer
the reader to \cite[p.384]{Kohlenbach}.\medskip

We now recall some mappings satisfying generalized contractive condition. A
mapping $T:X\rightarrow X$ is called:\newline
(i) Zamfirescu mapping \cite{Zamfirescu}\textbf{,} if there exist real
numbers $a,b$\ and $c$\ satisfying $a\in \left( 0,1\right) $ and $b,c\in
\left( 0,\frac{1}{2}\right) $\ such that for each pair of points $x,y$\ in $%
X,$ we have%
\begin{equation}
\begin{array}{l}
\text{(Z1) }d(Tx,Ty)\leq \alpha d(x,y) \\ 
\text{(Z2) }d(Tx,Ty)\leq b[d(x,Tx)+d(y,Ty)] \\ 
\text{(Z3) }d(Tx,Ty)\leq c[d(x,Ty)+d(y,Tx)];%
\end{array}
\label{2.1}
\end{equation}%
(ii) $C^{q}$-mapping \cite{Ciric}, if for some $h\in \lbrack 0,1)$ and for
all $x,y\in X,$ we have%
\begin{equation}
d(Tx,Ty)\leq h\max \{d(x,y),d(x,Tx),d(y,Ty),d(x,Ty),d(y,Tx)\};  \label{2.2}
\end{equation}%
(iii) generalized contractive mapping \cite{Naimpally}, if for some $h\in
\lbrack 0,1)$ and for all $x,y\in X,$ we have%
\begin{equation}
d(Tx,Ty)\leq h\max \{d(x,y),d(x,Tx),d(y,Ty),d(x,Ty)+d(y,Tx)\}.  \label{2.3}
\end{equation}%
(iv) generalized $C^{q}$-mapping \cite{Akhter}, if for some $h\in \lbrack
0,1)$ and for all $x,y\in X,$ we have%
\begin{equation}
d(Tx,Ty)\leq h\max \{d(x,y),d(x,Tx)+d(y,Ty),d(x,Ty)+d(y,Tx)\}.  \label{2.4}
\end{equation}%
\textbf{Remark 2.1.} It is evident from the above definitions that the class
of mappings defined in (2.4) contains properly the corresponding classes of
mappings defined in (2.1)-(2.3). However, the class of Zamfirescu mapping is
one of the most studied class of contractive mappings. For more on
contractive type mapping, we refer the reader to \cite{Browder}.\medskip

We now introduce different iterative algorithm, required in the sequel, in
Kohlenbach hyperbolic spaces. Let $T:X\rightarrow X$ be a given mapping and $%
x_{0}\in X$ be chosen arbitrarily, then the Picard, Mann, Ishikawa and
Xu-Noor iterative algorithms be defined, respectively, as follows:%
\begin{equation}
x_{n+1}=Tx_{n},  \label{2.5}
\end{equation}%
\begin{equation}
x_{n+1}=W(Tx_{n},x_{n},\alpha _{n}),  \label{2.6}
\end{equation}%
where $\{\alpha _{n}\}\subset \left[ 0,1\right] ,$%
\begin{eqnarray}
x_{n+1} &=&W(Ty_{n},x_{n},\alpha _{n})  \notag \\
y_{n} &=&W(Tx_{n},x_{n},\beta _{n}),  \label{2.7}
\end{eqnarray}%
where $\{\alpha _{n}\},$ $\{\beta _{n}\}\subset \left[ 0,1\right] ,$%
\begin{eqnarray}
x_{n+1} &=&W(Ty_{n},x_{n},\alpha _{n})  \notag \\
y_{n} &=&W(Tz_{n},x_{n},\beta _{n})  \label{2.8} \\
z_{n} &=&W(Tx_{n},x_{n},\gamma _{n}),  \notag
\end{eqnarray}%
where $\{\alpha _{n}\},$ $\{\beta _{n}\},\{\gamma _{n}\}\subset \left[ 0,1%
\right] .\medskip $

We now recall the concept introduced by Berinde \cite{Berinde1} for a
comparison of the rates of convergence of different iterative algorithms
involving a nonlinear mapping.

Let $\{a_{n}\}_{n=0}^{\infty },\{b_{n}\}_{n=0}^{\infty }$ be two sequences
of positive numbers that converge to $a,b,$ respectively. Assume that the
limit%
\begin{equation*}
l=\lim_{n\rightarrow \infty }\frac{\left\vert a_{n}-a\right\vert }{%
\left\vert b_{n}-b\right\vert },
\end{equation*}%
exists. If $l=0,$ then the sequence $\{a_{n}\}_{n=0}^{\infty }$ converges to 
$a$ faster than $\{b_{n}\}_{n=0}^{\infty }$ to $b.$ If $0<l<\infty ,$ then
we say that the two sequence $\{a_{n}\}_{n=0}^{\infty }$ and $%
\{b_{n}\}_{n=0}^{\infty }$ have the same rate of convergence. It is remarked
that the results concerning rates of convergence associated with the classes
of mappings defined in (2.1)-(2.3) have been established in \cite%
{Berinde,Berinde2,Berinde3,Rhoades}. See, also, \cite{Fukhar-Berinde} and
the references cited therein. We are now in a position to prove our main
results.

\section{Main Results}

This section is devoted to establish the results concerning iterative
construction of fixed points of the class of generalized $C^{q}$-mappings
and consequent rates of convergence for the modified Mann, Ishikawa and
Xu-Noor iterative algorithms in comparison to the classical Picard iterative
algorithm in Kohlenbach hyperbolic spaces.\newline
\textbf{Theorem 3.1.} Let $K$ be a nonempty closed convex subset of a
uniformly convex Kohlenbach hyperbolic space $X$ and let $T:K\rightarrow K$
be a generalized $C^{q}$-mapping. Assume that $F(T),$ the set of fixed
points of $T,$ is nonempty and the sequence $\left\{ \alpha _{n}\right\}
_{n=0}^{\infty }$ satisfies the following conditions:\newline
(C1): $0\leq \alpha _{n}<1$;\newline
(C2): $\tsum\limits_{n=0}^{\infty }\alpha _{n}=\infty .$\newline
Then the iterative algorithms defined in (2.5) and (2.6) converges to a
fixed point $p$ of $T$ provided that the iterative algorithms have same
initial guess $x_{0}\in K.$ Moreover, iterative algorithm defined in (2.5)
converges faster than (2.6) to the fixed point of $T.$\newline
\textbf{Proof.} Since $T$ is a generalized $C^{q}$-mapping, therefore, if $%
d(Tx,Ty)\leq h\{d(x,Tx)+d(y,Ty)\}$, then (2.4) becomes%
\begin{equation*}
d(Tx,Ty)\leq h\{d(x,Tx)+d(y,x)+d(x,Tx)+d(Tx,Ty)\}.
\end{equation*}%
So, we have%
\begin{equation}
d(Tx,Ty)\leq \frac{h}{1-h}\left\{ d(x,y)+2d(x,Tx)\right\} .  \label{3.1}
\end{equation}%
If $d(Tx,Ty)\leq h\{d(x,Ty)+d(y,Tx)\}$, then (2.4) becomes%
\begin{equation*}
d(Tx,Ty)\leq h\{d(x,Tx)+d(Tx,Ty)+d(y,x)+d(x,Tx)\}.
\end{equation*}%
Again, we have%
\begin{equation*}
d(Tx,Ty)\leq \frac{h}{1-h}\left\{ d(x,y)+2d(x,Tx)\right\} .
\end{equation*}%
Letting $\lambda =\max \left\{ h,\frac{h}{1-h}\right\} ,$ the above estimate
implies that%
\begin{equation}
d(Tx,Ty)\leq \lambda d(x,y)+2\lambda d(x,Tx).  \label{3.2}
\end{equation}%
Similarly, we can calculate the following inequality%
\begin{equation}
d(Tx,Ty)\leq \lambda d(x,y)+2\lambda d(y,Tx).  \label{3.3}
\end{equation}%
Let $p\in F(T),$ then it follows from the estimate (3.2) and the sequence
(2.5) that%
\begin{equation*}
d(x_{n+1},p)=d(Tx_{n},p)\leq \lambda d(x_{n},p).
\end{equation*}%
Continuing in this fashion, we have%
\begin{equation}
d(x_{n+1},p)\leq \lambda ^{n}d(x_{0},p).  \label{3.4}
\end{equation}%
Since $\lambda \in \lbrack 0,1),$ therefore, (3.4) implies that%
\begin{equation}
\lim_{n\rightarrow \infty }d(x_{n+1},p)=0.  \label{3.5}
\end{equation}%
Now utilizing the estimate (3.2) for the sequence (2.6), we get%
\begin{eqnarray}
d(x_{n+1},p) &\leq &\alpha _{n}d(Tx_{n},p)+(1-\alpha _{n})d(x_{n},p)  \notag
\\
&\leq &(1-(1-\lambda )\alpha _{n})d(x_{n},p).  \label{3.6}
\end{eqnarray}%
The estimate (3.6) inductively yields%
\begin{equation}
d(x_{n+1},p)\leq \prod\limits_{k=1}^{n}(1-(1-\lambda )\alpha _{k})d(x_{0},p).
\label{3.7}
\end{equation}%
Making use of conditions (C1) and (C2), the estimate (3.7) implies that%
\begin{equation}
\lim_{n\rightarrow \infty }d(x_{n+1},p)=0.  \label{3.8}
\end{equation}%
Hence the convergence of iterative algorithms (2.5) and (2.6) follows from
the estimates (3.5) and (3.8), respectively. In order to compare the rates
of convergence of iterative algorithms (2.5) and (2.6), we let $%
a_{n}=\lambda ^{n}$ and $b_{n}=\prod\limits_{k=1}^{n}(1-(1-\lambda )\alpha
_{k})d(x_{0},p).$\newline
Now, consider%
\begin{eqnarray*}
(1-\lambda )\alpha _{n} &\leq &(1-\lambda ) \\
-(1-\lambda ) &\leq &-(1-\lambda )\alpha _{n} \\
1-(1-\lambda ) &\leq &1-(1-\lambda )\alpha _{n} \\
\lambda &\leq &1-(1-\lambda )\alpha _{n} \\
\frac{\lambda }{1-(1-\lambda )\alpha _{n}} &\leq &1.
\end{eqnarray*}%
Moreover%
\begin{equation*}
\frac{\min \lambda }{\max [1-(1-\lambda )\alpha _{n}]}<1.
\end{equation*}%
Since $\prod\limits_{k=1}^{n}\frac{\lambda ^{k}}{[1-(1-\lambda )\alpha _{k}]}%
<\left( \frac{\min \lambda ^{k}}{\max [1-(1-\lambda )\alpha _{k}]}\right)
^{n},$ then%
\begin{equation*}
\lim_{n\rightarrow \infty }\frac{a_{n}}{b_{n}}=0.
\end{equation*}%
Hence (2.5) converges faster than (2.6) to the fixed point of $T.$\newline
\textbf{Theorem 3.2.} Let $K$ be a nonempty closed convex subset of a
uniformly convex Kohlenbach hyperbolic space $X$ and let $T:K\rightarrow K$
be a generalized $C^{q}$-mapping. Assume that $F(T),$ the set of fixed
points of $T,$ is nonempty and the sequences $\{\alpha _{n}\}_{n=0}^{\infty
} $ and $\{\beta _{n}\}_{n=0}^{\infty }$ satisfy the following conditions:%
\newline
(C1): $0\leq \alpha _{n},\beta _{n}<1$;\newline
(C2): $\tsum\limits_{n=0}^{\infty }\alpha _{n}=\infty .$\newline
Then the iterative algorithms defined in (2.6) and (2.7) converges to a
fixed point $p$ of $T$ provided that the iterative algorithms have same
initial guess $x_{0}\in K.$ Moreover, iterative algorithm defined in (2.6)
converges faster than (2.7) to the fixed point of $T.$\newline
\textbf{Proof.} Note that the convergence of (2.6) has already established
in Theorem 3.1. It remains to establish the convergence of (2.7) involving
the class of generalized $C^{q}$-mapping. For this, we proceed with the
following estimate:%
\begin{equation*}
d(x_{n+1},p)\leq \alpha _{n}d(Ty_{n},p)+(1-\alpha _{n})d(x_{n},p).
\end{equation*}%
On using (3.2), we get%
\begin{equation}
d(x_{n+1},p)\leq \alpha _{n}\lambda d(y_{n},p)+(1-\alpha _{n})d(x_{n},p).
\label{3.9}
\end{equation}%
Consider%
\begin{equation*}
d(y_{n},p)\leq \beta _{n}d(Tx_{n},p)+(1-\beta _{n})d(x_{n},p).
\end{equation*}%
Again, using (3.2), we get%
\begin{eqnarray*}
d(y_{n},p) &\leq &\beta _{n}\lambda d(x_{n},p)+(1-\beta _{n})d(x_{n},p) \\
&=&[\beta _{n}\lambda +(1-\beta _{n})]d(x_{n},p).
\end{eqnarray*}%
Substituting the above estimate in (3.9), we have%
\begin{eqnarray}
d(x_{n+1},p) &\leq &\alpha _{n}\lambda \lbrack \beta _{n}\lambda +(1-\beta
_{n})]d(x_{n},p)+(1-\alpha _{n})d(x_{n},p)  \notag \\
&=&[\alpha _{n}\beta _{n}\lambda ^{2}+\alpha _{n}\lambda (1-\beta
_{n})+(1-\alpha _{n})]d(x_{n},p)  \notag \\
&=&[1-\alpha _{n}(1-\lambda +\beta _{n}\lambda -\beta _{n}\lambda
^{2})]d(x_{n},p)  \notag \\
&=&[1-\alpha _{n}((1-\lambda )+\beta _{n}\lambda (1-\lambda ))]d(x_{n},p) 
\notag \\
&=&[1-\alpha _{n}(1-\lambda )(1+\beta _{n}\lambda )]d(x_{n},p).  \label{3.10}
\end{eqnarray}%
Consider%
\begin{eqnarray*}
1-\lambda &\leq &1+\beta _{n}\lambda \\
\alpha _{n}(1-\lambda )(1-\lambda ) &\leq &\alpha _{n}(1-\lambda )(1+\beta
_{n}\lambda ) \\
-\alpha _{n}(1-\lambda )(1+\beta _{n}\lambda ) &\leq &-\alpha _{n}(1-\lambda
)(1-\lambda ) \\
1-\alpha _{n}(1-\lambda )(1+\beta _{n}\lambda ) &\leq &1-\alpha
_{n}(1-\lambda )^{2}.
\end{eqnarray*}%
Utilizing the above assertion, the estimate (3.10) implies that%
\begin{equation}
d(x_{n+1},p)\leq \lbrack 1-\alpha _{n}(1-\lambda )^{2}]d(x_{n},p).
\label{3.11}
\end{equation}%
Continuing in this fashion, we have%
\begin{equation*}
d(x_{n+1},p)\leq \tprod\limits_{k=1}^{n}[1-\alpha _{k}(1-\lambda
)^{2}]d(x_{0},p).
\end{equation*}%
Using the fact that $\lambda \in \lbrack 0,1)$ and conditions (C1)-(C2), we
get%
\begin{equation}
\lim_{n\rightarrow \infty }d(x_{n+1},p)=0.  \label{3.12}
\end{equation}%
The estimate (3.12) implies that (2.7) converges the fixed point $p$ of $T.$
In order to compare the rates of convergence of (2.6) and (2.7), we must
compare $a_{n}=\tprod\limits_{k=1}^{n}[1-\alpha _{k}(1-\lambda )]$ and $%
b_{n}=\tprod\limits_{k=1}^{n}[1-\alpha _{k}(1-\lambda )^{2}].$ For this, we
reason as follow:%
\begin{eqnarray*}
\alpha _{k}(1-\lambda )(1-\lambda ) &\leq &\alpha _{k}(1-\lambda ) \\
-\alpha _{k}(1-\lambda ) &\leq &-\alpha _{k}(1-\lambda )^{2} \\
1-\alpha _{k}(1-\lambda ) &\leq &1-\alpha _{k}(1-\lambda )^{2} \\
\frac{1-\alpha _{k}(1-\lambda )}{1-\alpha _{k}(1-\lambda )^{2}} &\leq &1.
\end{eqnarray*}%
Also%
\begin{equation*}
\frac{\min \{1-\alpha _{k}(1-\lambda )\}}{\max \{1-\alpha _{k}(1-\lambda
)^{2}\}}<1.
\end{equation*}%
Since $\dprod\limits_{k=1}^{n}\frac{[1-\alpha _{k}(1-\lambda )]}{[1-\alpha
_{k}(1-\lambda )^{2}]}<$ $\left( \frac{\min \{1-\alpha _{k}(1-\lambda )\}}{%
\max \{1-\alpha _{k}(1-\lambda )^{2}\}}\right) ^{n},$ then%
\begin{equation*}
\lim_{n\rightarrow \infty }\frac{a_{n}}{b_{n}}=0.
\end{equation*}%
Hence (2.6) converges faster than (2.7) to the fixed point of $T.$\newline
\textbf{Theorem 3.3.} Let $K$ be a nonempty closed convex subset of a
uniformly convex Kohlenbach hyperbolic space $X$ and let $T:K\rightarrow K$
be a generalized $C^{q}$-mapping. Assume that $F(T),$ the set of fixed
points of $T,$ is nonempty and the sequences $\{\alpha _{n}\}_{n=0}^{\infty
},\{\beta _{n}\}_{n=0}^{\infty }$ and $\{\gamma _{n}\}_{n=0}^{\infty }$
satisfy the following conditions:\newline
(C1): $0\leq \alpha _{n},\beta _{n},\gamma _{n}<1$;\newline
(C2): $\tsum\limits_{n=0}^{\infty }\alpha _{n}=\infty .$\newline
Then the iterative algorithms defined in (2.7) and (2.8) converges to a
fixed point $p$ of $T$ provided that the iterative algorithms have same
initial guess $x_{0}\in K.$ Moreover, iterative algorithm defined in (2.7)
converges faster than (2.8) to the fixed point of $T.$\newline
\textbf{Proof.} Note that the convergence of (2.7) has already established
in Theorem 3.2. It remains to establish the convergence of (2.8) involving
the class of generalized $C^{q}$-mapping. For this, we proceed with the
following estimates:%
\begin{eqnarray}
d(z_{n},p) &=&d(W(Tx_{n},x_{n},\gamma _{n}),p)  \notag \\
&\leq &\gamma _{n}d(p,Tx_{n})+(1-\gamma _{n})d(x_{n},p)  \notag \\
&\leq &(1-\gamma _{n}(1-\lambda ))d(x_{n},p)  \label{3.13}
\end{eqnarray}%
and%
\begin{eqnarray}
d(y_{n},p) &=&d(W(Tz_{n},x_{n},\beta _{n}),p)  \notag \\
&\leq &\beta _{n}d(p,Tz_{n})+(1-\beta _{n})d(x_{n},p)  \notag \\
&\leq &\beta _{n}\lambda d(z_{n},p)+(1-\beta _{n})d(x_{n},p).  \label{3.14}
\end{eqnarray}%
Substituting (3.13) in (3.14), we have%
\begin{equation}
d(y_{n},p)\leq \beta _{n}\lambda \left[ (1-\gamma _{n}(1-\lambda ))d(x_{n},p)%
\right] +(1-\beta _{n})d(x_{n},p)  \label{3.15}
\end{equation}%
Moreover%
\begin{eqnarray}
d(x_{n+1},p) &=&d(W(Ty_{n},x_{n},\alpha _{n}),p)  \notag \\
&\leq &\alpha _{n}d(p,Ty_{n})+(1-\alpha _{n})d(x_{n},p)  \notag \\
&\leq &\alpha _{n}\lambda d(y_{n},p)+(1-\alpha _{n})d(x_{n},p)  \label{3.16}
\end{eqnarray}%
Substituting (3.15) in (3.16),\ we get%
\begin{eqnarray*}
d(x_{n+1},p) &\leq &\left[ \alpha _{n}\lambda \left\{ \beta _{n}\lambda
(1-\gamma _{n}(1-\lambda ))+1-\beta _{n}\right\} +(1-\alpha _{n})\right]
d(x_{n},p) \\
&=&\left\{ 1-\alpha _{n}\left[ 1-\beta _{n}\lambda ^{2}(1-\gamma _{n}+\gamma
_{n}\lambda )-\lambda (1-\beta _{n})\right] \right\} d(x_{n},p) \\
&=&\left\{ 1-\alpha _{n}\left[ 1-\beta _{n}\lambda ^{2}+\beta _{n}\gamma
_{n}\lambda ^{2}-\beta _{n}\gamma _{n}\lambda ^{3}-\lambda +\beta
_{n}\lambda \right] \right\} d(x_{n},p) \\
&=&\left\{ 1-\alpha _{n}\left[ 1-\lambda +(1-\lambda )(\beta _{n}\lambda
+\beta _{n}\gamma _{n}\lambda ^{2})\right] \right\} d(x_{n},p) \\
&=&\left\{ 1-\alpha _{n}(1-\lambda )\left[ 1+\beta _{n}\gamma _{n}\lambda
^{2}+\beta _{n}\lambda \right] \right\} d(x_{n},p) \\
&\leq &\left\{ 1-\alpha _{n}(1-\lambda )\right\} d(x_{n},p).
\end{eqnarray*}%
Making use of conditions (C1) and (C2), the above estimate implies that%
\begin{equation}
\lim d(x_{n+1},p)=0.  \label{3.17}
\end{equation}%
Now we use the estimate (3.3) for the iterative algorithm (2.8) to get the
following estimates:%
\begin{eqnarray}
d(z_{n},p) &=&d(W(Tx_{n},x_{n},\gamma _{n}),p)  \notag \\
&\leq &\gamma _{n}d(Tx_{n},p)+(1-\gamma _{n})d(x_{n},p)  \notag \\
&\leq &3\lambda \gamma _{n}d(x_{n},p)+(1-\gamma _{n})d(x_{n},p)  \notag \\
&=&\left[ 3\lambda \gamma _{n}+(1-\gamma _{n})\right] d(x_{n},p),
\label{3.18}
\end{eqnarray}%
and%
\begin{equation}
d(y_{n},p)\leq 3\lambda \beta _{n}d(z_{n},p)+(1-\beta _{n})d(x_{n},p).
\label{3.19}
\end{equation}%
Substituting (3.18) in (3.19), we get%
\begin{equation}
d(y_{n},p)\leq 3\lambda \beta _{n}\left[ 3\lambda \gamma _{n}+(1-\gamma
_{n})+(1-\beta _{n})\right] d(x_{n},p).  \label{3.20}
\end{equation}%
Now, consider%
\begin{equation}
d(x_{n+1},p)\leq 3\lambda \alpha _{n}d(y_{n},p)+(1-\alpha _{n})d(x_{n},p).
\label{3.21}
\end{equation}%
Substituting (3.20)\ in (3.21) and then simplifying the terms, we have%
\begin{eqnarray}
d(x_{n+1},p) &\leq &\left[ 1-\alpha _{n}\left( 1-3\lambda \right) \left\{
1+9\lambda ^{2}\beta _{n}\gamma _{n}+3\lambda \beta _{n}\right\} \right]
d(x_{n},p)  \notag \\
&\leq &\left[ 1-\alpha _{n}\left( 1-3\lambda \right) \right] d(x_{n},p).
\label{3.22}
\end{eqnarray}%
Again, making use of conditions (C1) and (C2), the above estimate implies
that%
\begin{equation}
\lim d(x_{n+1},p)=0.  \label{3.23}
\end{equation}%
In order to compare the rates of convergence of (2.7) and (2.8), we must
compare $a_{n}=\tprod\limits_{k=1}^{n}[1-\alpha _{k}(1-\lambda )^{2}]$ and $%
b_{n}=\tprod\limits_{k=1}^{n}[1-\alpha _{k}(1-3\lambda )].$ For this, we
have the following two cases:\newline
Case (I). Let $\lambda \in \lbrack 0,\frac{1}{3}],$ then $a_{n}\leq 1$ and $%
b_{n}=1,$ therefore, we have $\lim_{n\rightarrow \infty }(\frac{a_{n}}{b_{n}}%
)=0$.\newline
Case (II). Let $\lambda \in (\frac{1}{3},1),$ then again $a_{n}\leq 1$ and%
\begin{equation*}
b_{n}=\dprod\limits_{k=1}^{n}\left[ 1-\alpha _{k}\left( 1-3\lambda \right)
\left\{ 1+9\lambda ^{2}\beta _{k}\gamma _{k}+3\lambda \beta _{k}\right\} %
\right] \geq 1.
\end{equation*}%
So%
\begin{equation*}
\frac{a_{n}}{b_{n}}=\dprod\limits_{k=1}^{n}\left[ \frac{1-\alpha
_{k}(1-\lambda )^{2}}{1-\alpha _{k}\left( 1-3\lambda \right) \left\{
1+9\lambda ^{2}\beta _{k}\gamma _{k}+3\lambda \beta _{k}\right\} }\right]
\leq 1.
\end{equation*}%
Consequently%
\begin{equation*}
\frac{\min \{1-\alpha _{k}(1-\lambda )^{2}\}}{\max \{1-\alpha _{k}\left(
1-3\lambda \right) \left\{ 1+9\lambda ^{2}\beta _{k}\gamma _{k}+3\lambda
\beta _{k}\right\} \}}<1.
\end{equation*}%
Since $\dprod\limits_{k=1}^{n}\left[ \frac{1-\alpha _{k}(1-\lambda )^{2}}{%
1-\alpha _{k}\left( 1-3\lambda \right) \left\{ 1+9\lambda ^{2}\beta
_{k}\gamma _{k}+3\lambda \beta _{k}\right\} }\right] <\left( \frac{\min
\{1-\alpha _{k}(1-\lambda )^{2}\}}{\max \{1-\alpha _{k}\left( 1-3\lambda
\right) \left\{ 1+9\lambda ^{2}\beta _{k}\gamma _{k}+3\lambda \beta
_{k}\right\} \}}\right) ^{n},$ therefore, we get%
\begin{equation*}
\lim_{n\rightarrow \infty }\frac{a_{n}}{b_{n}}=0.
\end{equation*}%
This implies that, in both cases, (2.7) converges faster than (2.8) to the
fixed point of $T.$\newline
\textbf{Remark 3.4.} As an applications of Theorems (3.1)-(3.3), we can
establish similar kind of results for the classes of mappings defined in
(2.1)-(2.3) in Kohlenbach hyperbolic spaces. As a consequence, our results
generalize the corresponding results from linear spaces to more general
setup of spaces.

\end{document}